\documentclass[a4paper,11pt]{delfim}

\begin{document}

%%%%%%%%%%%%%%%%%%%%%%%%%%%%%%%%%%%%%%%%%%%%%%%%%%%%%%

\title{Regularity of Solutions to Second-Order Integral
Functionals in Variational Calculus}

\author{\textbf{Moulay Rchid Sidi Ammi$^1$ and Delfim F. M. Torres$^2$}}

\date{$^1$Department of Mathematics\\
University of Aveiro\\
3810-193 Aveiro, Portugal\\
sidiammi@mat.ua.pt\\[0.3cm]
$^2$Department of Mathematics\\
University of Aveiro\\
3810-193 Aveiro, Portugal\\
delfim@ua.pt}

\maketitle

%%%%%%%%%%%%%%%%%%%%%%%%%%%%%%%%%%%%%%%%%%%%%%%%%%%%%%

\begin{abstract}
\noindent \emph{We obtain regularity conditions of a new type of problems of the
calculus of variations with second-order derivatives. As a
corollary, we get non-occurrence of the Lavrentiev phenomenon. Our
main result asserts that autonomous integral functionals of the
calculus of variations with a Lagrangian having superlinearity partial
derivatives with respect to the higher-order derivatives admit only
minimizers with essentially bounded derivatives.
}

\medskip


\medskip


\end{abstract}

%%%%%%%%%%%%%%%%%%%%%%%%%%%%%%%%%%%%%%%%%%%%%%%%%%%%%%%%%%%%%%%%%%%%

\section{Introduction}

Let $\mathcal{L}(t, x^{0},\ldots, x^{m})$ be a given $C^{1}([a, b]
\times \mathbb{R}^{(m+1)\times n})$ real valued function. The
problem of the calculus of variations with high-order derivatives
consists in minimizing an integral functional
\begin{equation}
\label{Pm} J^{m}[x(\cdot)]= \int_{a}^{b}\mathcal{L}\left(t, x(t), \dot{x}(t),
\ldots, x^{(m)}(t)\right) dt  \tag{$P_{m}$}
\end{equation}
over a certain class $\mathcal{X}$ of functions $x : [a, b]
\rightarrow \mathbb{R}^{n}$ satisfying the boundary conditions
\begin{equation}
\label{eq:boundaryCond} x(a)= x_{a}^{0}\, ,  x(b)= x_{b}^{0}\, ,
\ldots\, , x^{(m-1)}(a)= x_{a}^{m-1}\, , x^{(m-1)}(b)= x_{b}^{m-1}
\, .
\end{equation}
Often it is convenient to write $x^{(1)}= x'$, $x^{(2)}= x''$, and
sometimes we revert to the standard notation used in mechanics:
$x'=\dot{x}$, $x''= \ddot{x}$. Such problems arise, for instance, in
connection with the theory of beams and rods \cite{smi}. Further,
many problems in the calculus of variations with higher-order
derivatives describe important optimal control problems with linear
dynamics \cite{SarychevTorres2}.

Regularity theory for optimal control problems is a fertile field of
research and a source of many challenging mathematical issues and
interesting applications \cite{cla,torres3,torres2}. The essential
points in the theory are: (i) existence of minimizers and (ii)
necessary optimality conditions to identify those minimizers.

The first systematic approach to existence theory was introduced
by Tonelli in 1915 \cite{ton}, who showed that existence of
minimizers is guaranteed in the Sobolev space $W^{m}_{m}$ of
absolutely continuous functions. The direct method of Tonelli
proceeds in three steps: (i) regularity and convexity with respect
to the highest-derivative of the Lagrangian $\mathcal{L}$ guarantees lower
semi-continuity, (ii) the coercivity condition (the Lagrangian $\mathcal{L}$
must grow faster than a linear function) insure compactness, (iii)
by the compactness principle, one gets the existence of minimizers
for the problem \eqref{Pm}. Typically, Tonelli's existence theorem
for \eqref{Pm} is formulated as follows
\cite{cla,clavin90}: under hypotheses (H1)-(H3) on the Lagrangian
$\mathcal{L}$,
\begin{enumerate}
\item[(H1)] $\mathcal{L}(t, x^{0},\ldots, x^{m})$ is locally Lipschitz in $(t, x^{0},\ldots, x^{m})$;
\item[(H2)] $\mathcal{L}(t, x^{0},\ldots, x^{m})$ is convex as a function of the last argument $x^{m}$;
\item[(H3)] $\mathcal{L}(t, x^{0},\ldots, x^{m})$ is coercive in $x^{m}$, \textrm{i.e.}
$\exists$ $\Theta : [0, \infty)\rightarrow \mathbb{R}$ such that
\begin{gather*}
\lim_{r \rightarrow \infty} \frac{\Theta (r)}{r}=+\infty \, , \\
\mathcal{L}(t, x^{0},\ldots, x^{m}) \geq \Theta (|x^{m}|) \mbox{ for all } (t,
x^{0},\ldots, x^{m}) \, ,
\end{gather*}
\end{enumerate}
there exists a minimizer to problem \eqref{Pm} in the class
$W^{m}_{m}$.

The main necessary condition in optimal control is the famous
Pontryagin maximum principle, which includes all the classical
necessary optimality conditions of the calculus of variations
\cite{Pontryagin}. It turns out that the hypotheses (H1)-(H3) do not
assure the applicability of the necessary optimality conditions,
being required more regularity on the class of admissible functions
\cite{BallMizel}. For \eqref{Pm}, the Pontryagin maximum principle
\cite{Pontryagin} is established assuming $x \in W_{m}^{\infty}
\subset W^{m}_{m}$.

In the case $m=1$, extra information about the minimizers was
proved, for the first time, by Tonelli himself \cite{ton}. Tonelli
established that, under the hypotheses (H2) and (H3) of convexity
and coercivity, the minimizers $x$ have the property that $\dot{x}$
is locally essentially bounded on an open subset $\Omega \subset [a,
b]$ of full measure. If
\begin{equation}
\label{eq:conTonelliMorrey} \left|\frac{\partial \mathcal{L}}{\partial
x}\right| + \left|\frac{\partial \mathcal{L}}{\partial \dot{x}}\right| \le c
|\mathcal{L}| + r \, ,
\end{equation}
for some constants $c$ and $r$, $c > 0$, then $\Omega = [a,b]$
($\dot{x}(t)$ is essentially bounded in all points $t$ of $[a,b]$,
\textrm{i.e.} $x \in W_{1}^{\infty}$), and the Pontryagin maximum
principle, or the necessary condition of Euler-Lagrange, hold.
Condition \eqref{eq:conTonelliMorrey} is now known in the literature
as the Tonelli-Morrey regularity condition
\cite{FrancisClarke,clavin,SarychevTorres2}. Since L.~Tonelli and
C.~B.~Morrey, several Lipschitzian regularity conditions were
obtained for the problem \eqref{Pm} with $m = 1$: S.~Bernstein (for
the scalar case $n = 1$), F.~H.~Clarke and R.~B.~Vinter (for the
vectorial case $n > 1$) obtained \cite{clavin85} the condition
\begin{equation*}
\left|\left(\frac{\partial^2 \mathcal{L}}{\partial \dot{x}^2}\right)^{-1}
\left( \frac{\partial \mathcal{L}}{\partial x} -\frac{\partial^2 \mathcal{L}}{\partial
\dot{x} \partial t} -\frac{\partial^2 \mathcal{L}}{\partial \dot{x} \partial
x} \dot{x}\right)\right| \le c \left(\left|\dot{x}\right|^3 +
1\right) \, , \quad \frac{\partial^2 \mathcal{L}}{\partial \dot{x}^2} > 0 \, ;
\end{equation*}
F.~H.~Clarke and R.~B.~Vinter \cite{clavin85} the regularity
conditions
\begin{equation}
\label{eq:aut} \left|\frac{\partial \mathcal{L}}{\partial t}\right| \le c
\left|\mathcal{L}\right| + k(t) \, , \quad k(\cdot) \in L_1 \, ,
\end{equation}
and
\begin{equation*}
\left|\frac{\partial \mathcal{L}}{\partial x}\right| \le c \left|\mathcal{L}\right| +
k(t) \left|\frac{\partial \mathcal{L}}{\partial \dot{x}}\right| + m(t) \, , \quad
k(\cdot), \, m(\cdot) \in L_1 \, ;
\end{equation*}
and A.~V.~Sarychev and D.~F.~M.~Torres \cite{SarychevTorres1} the condition
\begin{equation}
\label{eq:MSc} \left(\left|\frac{\partial \mathcal{L}}{\partial t}\right| +
\left|\frac{\partial \mathcal{L}}{\partial x}\right| \right)
\left|\dot{x}\right|^{\mu} \le \gamma \mathcal{L}^\beta + \eta \, , \quad
\gamma > 0 \, , \beta < 2 \, , \mu \ge \max
\left\{\beta-1,-1\right\} \, .
\end{equation}
Lipschitzian regularity theory for the problem of the calculus of
variations with $m = 1$ is now a vast discipline (see \textrm{e.g.}
\cite{Cellina,CellinaFerriero,DalMaso,Ornelas,torres3} and
references therein). Results for $m > 1$ are scarcer: we are
aware of the results in \cite{clavin90,SarychevTorres1,torres1}. In
1997 A.V.~Sarychev \cite{Sarychev} proved that the second-order
problems of the calculus of variations may show new phenomena
non-present in the first-order case: under the hypotheses (H1)-(H3) of
Tonelli's existence theory, autonomous problems \eqref{Pm} with $m =
2$ may present the Lavrentiev phenomenon \cite{Lav}. This is not a
possibility for $m = 1$, as shown by the Lipschitzian regularity
condition \eqref{eq:aut}. Sarychev's result was recently extended by
A.~Ferriero \cite{Ferriero} for the case $m > 2$. It is also shown
in \cite{Ferriero} that, under some standard hypotheses, the
problems of the calculus of variations \eqref{Pm} with Lagrangians
only depending on two consecutive derivatives $x^{(\gamma)}$ and
$x^{(\gamma+1)}$, $\gamma \ge 0$, do not exhibit the Lavrentiev
phenomenon for any boundary conditions \eqref{eq:boundaryCond} (for
$m = 1$, this follows immediately from \eqref{eq:aut}). In the case
in which the Lagrangian only depends on the higher-order derivative
$x^{(m)}$, it is possible to prove more
\cite[Corollary~2]{SarychevTorres1}: when $\mathcal{L} =
\mathcal{L}\left(x^{(m)}\right)$, all the minimizers predicted by the
existence theory belong to the space $W_{m}^{\infty} \subset
W^{m}_{m}$ and satisfy the Pontryagin maximum principle
(regularity). As to whether this is the case or not for Ferriero's
problem with Lagrangians only depending on consecutive derivatives
$x^{(\gamma)}$ and $x^{(\gamma+1)}$, seems to be an open question.

The results of Sarychev \cite{Sarychev} and Ferriero \cite{Ferriero}
on the Lavrentiev phenomenon show that the problems of the calculus
of variations with higher-order derivatives are richer than the
problems with $m = 1$, but also show, in our opinion, that the
regularity theory for higher-order problems is underdeveloped. One
can say that the Lipschitzian regularity conditions found in the
literature for the higher-order problems of the calculus of
variations are a generalization of the above mentioned conditions
for $m = 1$: \cite{clavin90} generalizes \eqref{eq:conTonelliMorrey}
for $m > 1$, \cite{SarychevTorres1} generalizes \eqref{eq:MSc} for problems
of optimal control with control-affine dynamics, \cite{torres1}
generalizes \eqref{eq:conTonelliMorrey} for optimal control problems
with more general nonlinear dynamics. To the best of our knowledge,
there exist no regularity conditions for the higher-order problems
of the calculus of variations of a different type from those also
obtained (also valid) for the first-order problems. We give here
what we claim to be a new regularity condition which is of a
different nature than those appearing for the first-order problems.
For the sake of simplicity, we restrict ourselves to second-order
problems ($m = 2$). The results of the paper can be naturally
extended to derivatives of higher order than two, but the proofs
become rather technical. While existence follows by imposing
coercivity to the Lagrangian $\mathcal{L}$ (hypothesis (H3)), we
prove (\textrm{cf.} Theorem~\ref{theorem4}) that for the autonomous
second-order problems of the calculus of variations, regularity
follows by imposing a superlinearity condition to the partial
derivatives $\frac{\partial \mathcal{L}}{\partial \ddot{x}_i}$ of
the Lagrangian. We observe that our condition is intrinsic to the
higher-order problems: for autonomous problems of the calculus of
variations with $m = 1$ \eqref{eq:aut} is trivially satisfied and no
superlinearity on the partial derivatives of $\mathcal{L}$ is needed,
while such conditions are required in the higher-order case as a
consequence of Sarychev's results \cite{Sarychev}.

%%%%%%%%%%%%%%%%%%%%%%%%%%%%%%%%%%%%%%%%%%%%

\section{Outline of the paper and hypotheses}

We shall limit ourselves to the second order problems of the calculus of
variations, \textrm{i.e.} to the problem of minimizing
\begin{equation}
\label{P2} \int_{a}^{b} \mathcal{L}\left(t, x(t),\dot{x}(t) ,
\ddot{x}(t)\right) dt  \tag{$P_{2}$}
\end{equation}
for some given Lagrangian $\mathcal{L}(\cdot,\cdot,\cdot,\cdot)$, assumed to
be a $C^1$ function with respect to all arguments. In this case it
is appropriate to choose the admissible functions $x$ to be twice
continuously differentiable with derivatives $\dot{x}$ and
$\ddot{x}$ in $L^{2}$, i.e $\mathcal{X} = W^{2}_{2}$. In
Section~\ref{sec:nec:cond} we establish generalized integral forms
of duBois-Reymond and Euler-Lagrange necessary conditions valid for
$\mathcal{X} = W^{2}_{2}$ (the optimal solutions $x$ may have unbounded
derivatives $\dot{x}$, $\ddot{x}$). Then, in Section~\ref{sec:Reg},
we obtain regularity conditions under which all the minimizers of
\eqref{P2} are in $W_{2}^{\infty} \subset W^{2}_{2}$ and thus
satisfy the classical necessary conditions. In general terms, the
techniques used here are extensions of those appearing in \cite{ces}
for one-derivative problems.

\medskip

In the sequel we shall assume the following hypotheses:
\begin{description}

\item[$(S_{0})$] There exists a continuous function $S(t, s, v,
w)\geq 0$, $(t, s, v, w) \in \mathbb{R}^{1+3n}$, and some $\delta
> 0$, such that the function $t \rightarrow S(t, x(t), x'(t),
x''(t))$ is $L^{2}$-integrable in $[a, b]$ and
$$
\left|\frac{\partial \mathcal{L}}{\partial t}(\tau, x, x', x'')\right| \leq
S(t, x, x', x''),
$$
for all $t \in [a, b]$, $\left|\tau -t\right|< \delta, x=x(t)$.

\item[$(S_{i})$] There exists a nonnegative continuous function
$G(\cdot,\cdot,\cdot,\cdot)$, and some $\delta >0$, such that the
function $t \rightarrow G(t, x(t), x'(t), x''(t))$ is
$L^{2}$-integrable on $[a, b]$, and
\begin{equation*}
\begin{split}
\left|\frac{\partial \mathcal{L}}{\partial x_{i}}(t, y, x', x'')\right| &\leq G(t, x, x', x'')\, ,\\
\left|\frac{\partial \mathcal{L}}{\partial \dot{x}_{i}}(t, x, y, x'')\right| &\leq G(t, x, x', x'')\, ,\\
\left|\frac{\partial \mathcal{L}}{\partial \ddot{x}_{i}}(t, x, x',
y)\right| &\leq G(t, x, x', x'')\, ,
\end{split}
\end{equation*}
for all $t \in [a, b]$, $x$, $x'$, $x'' \in \mathbb{R}^{n}$,
$x=(x_{1},\ldots, x_{n})\in \mathbb{R}^{n}$,
$y=(y_{1},\ldots, y_{n})\in \mathbb{R}^{n}$, $y_{j}= x^{(k)}_{j}(t)$
for $j\neq i$,  $ \left|y_{i}-x_{i}^{(k)}(t)\right|\leq \delta$,
$i=1, \ldots, n$ and $k=0, 1, 2$, where $x_{i}^{(k)}(t)$ is the
$i^{th}$ component of the $k^{th}$ derivative with the convention
$x_{i}^{(0)}(t)=x_{i}(t)$.
\end{description}

\begin{remark}
Hypothesis $(S_{0})$ is certainly verified if $\mathcal{L}(t, x, \dot{x},
\ddot{x})$ does not depend on $t$: $(S_{0})$ holds trivially in the
autonomous case. Conditions $(S_{i})$, $i = 0,\ldots,n$, are needed
in the proof of Theorems~\ref{theorem1} and \ref{theorem2}
to justify the usual rule of differentiation under the sign of an integral.
\end{remark}

%%%%%%%%%%%%%%%%%%%%%

\section{Generalized integral form of duBois-Reymond and Euler-Lagrange equations}
\label{sec:nec:cond}

In this section we prove integral forms of duBois Reymond
and Euler-Lagrange equations (see \eqref{eq2} and \eqref{eq6} below,
respectively). For this, we consider an arbitrary change of the
independent variable $t$. Let $s$ be the arc length parameter on the
curve $C_{0}: x = x(t)$, $a \leq t  \leq b$, so that the Jordan
length of $C_{0}$ is $s(t)=\int_{a}^{t} \sqrt{
1+|x'(\tau)|^{2}}d\tau$ with $s(a)=0$, $s(b)=l$ and $s(t)$ is
absolutely continuous with $s'(t)\geq 1$ $a.e.$ Thus $s(t)$ and its
inverse $t(s)$, $0\leq s \leq l$, are absolutely continuous with
$t'(s)>0$ $a.e.$ in $[0, l]$. If $X(s)= x(t(s))$, $0 \leq s \leq l$,
then $t(s)$ and $X(s)$ are Lipschitzian of constant one in $[0, l]$.
By change of variable,
\begin{equation*}
\begin{split}
I[x] &= \int_{a}^{b}\mathcal{L}\left(t, x(t),\dot{x}(t) , \ddot{x}(t) \right) dt \\
&= \int_{0}^{l} \mathcal{L}\left(t(s), X(s),\frac{X'(s)}{t'(s)} ,
\frac{1}{t^{'2}(s)}\left(X''(s)- \frac{X'(s)}{t'(s)}t''(s)\right)\right)t'(s) ds
\, .
\end{split}
\end{equation*}
Setting $F(t, x, t', x', t'', x'') = \mathcal{L}\left(t,
x,\frac{x'}{t'},\frac{1}{t^{'2}}(x''- \frac{x'}{t'}t'')\right)t'$,
we have:
$$
I[x]= J[C]= J[X] =\int_{0}^{l}F\left(t(s), X(s), t'(s), X'(s),
t''(s), X''(s)\right) ds \, .
$$

%%%%%%%%%%%%%%%%%%%%%%%%%%%%%%%%%%%%%%%%%%%%

\subsection{Generalized duBois-Reymond equations}

The following necessary condition will be useful to prove our
regularity theorem (Theorem~\ref{theorem4} on
Section~\ref{sec:Reg}).

\begin{theorem}\label{theorem1}
Under hypotheses $(S_{i})_{0 \leq i \leq n}$,
if $x(\cdot) \in W_2^2$ is a minimizer of problem \eqref{P2}, then
the following integral
form of duBois-Reymond necessary condition holds:
\begin{equation}
\phi_{0}(s)= \frac{\partial F}{\partial t''}(\theta(s))
- \int_{0}^{s}\frac{\partial F}{\partial t'}(\theta(\sigma)) d\sigma +
\int_{0}^{s}\int_{0}^{\tau}\frac{\partial F}{\partial t}(\theta(\sigma))
d\sigma d\tau = c_{0},
\quad 0\leq \tau \leq s \leq l
 , \label{eq2}
\end{equation}
where functions $\frac{\partial F}{\partial t''}$, $\frac{\partial
F}{\partial t'}$, $\frac{\partial F}{\partial t}$ are evaluated at
$\theta(s) = (t(s), X(s), t'(s), X'(s), t''(s), X''(s))$
and $c_{0}$ is a constant.
\end{theorem}
\begin{proof}
It is to be noted that $(t(s), X(s), t'(s), X'(s), t''(s))$ may not
exist in a set of null-measure of all $s$. The proof is done by
contradiction. Suppose that \eqref{eq2} is not true. Then, there
exist constants $d_{1}<d_{2}$ and disjoints sets $E_{1}^{*}$ and
$E_{2}^{*}$ of non-zero measure such that
\begin{gather}
\phi_{0}(s)\leq d_{1} \mbox{ for } s \in
E_{1}^{*}, \nonumber\\
\phi_{0}(s)\geq  d_{2} \mbox{ for } s \in E_{2}^{*}, \nonumber
\end{gather}
while $t'(s)>0$ $a.e$ in $[0, l]$. Hence there exist some
constant $k>0$ and two subsets $E_{1}$, $E_{2}$ of positive measure
of $E_{1}^{*}$, $E_{2}^{*}$, such that
\begin{gather}
t'(s)\geq k>0, \quad \phi_{0}(s)\leq d_{1} \quad \mbox{ for } s \in
E_{1}, \quad \left|E_{1}\right|>0 \, , \label{eq3}\\
t'(s)\geq k>0,  \quad \phi_{0}(s)\geq  d_{2}  \quad \mbox{ for } s
\in E_{2},  \quad \left|E_{2}\right|>0 \, . \label{eq4}
\end{gather}
Let us consider
$$
\psi (s)= \int_{0}^{s}\int_{0}^{\tau} \{ \left|E_{2}\right| \chi_{1}(\sigma)
-\left|E_{1}\right| \chi_{2}(\sigma)\} d\sigma d\tau\, , \quad
0\leq \tau \leq  s \leq l \, ,
$$
where $\chi_{i}$ denotes the indicator function defined by
 $$
 \chi_{i}(s)=
 \left\{
 \begin{array}{rll}
  1  & \mbox{ for } & s \in E_{i},
\\
  0  & \mbox{ for } & s \in [0, l]/E_{i}, \quad  i=1, 2 \mbox{ and }  0\leq s \leq
  l.
   \end{array}
    \right.
   $$
We have that $\psi'$ is an absolutely continuous function in $[0,
l]$ with $\psi'(0)= \psi'(l)= 0$. Moreover,
$$
 \psi''(s)=
 \left\{
 \begin{array}{rll}
  -\left|E_{1}\right| & \mbox{ $a.e$ } & s \in E_{2} \, ,
\\
 \left|E_{2}\right| & \mbox{ $a.e$ } & s \in E_{1} \, , \\
0 & \mbox{ $a.e$ } & s \in [0, l]-E_{1} \bigcup E_{2} \, .
   \end{array}
    \right.
   $$
We also define $C_{\alpha}: t= t_{\alpha}(s)$, $x= X_{\alpha}(s)$,
$0 \leq s  \leq l$, by setting
$$
t_{\alpha}(s)= t(s)+ \alpha \psi (s) +\alpha^{2} \psi' (s) \, ,
$$
$$
X_{\alpha}(s)= X(s), \quad 0 \leq s \leq l, \quad
\left|\alpha\right|\leq 1 \, .
$$
Let $\rho >0$ be chosen in such a way that  $t, \tau \in [a, b]$ and
$\left|t-\tau\right|< \rho$ imply $|x(t)-x(\tau)|\leq \delta$,
where $\delta$ is the constant in condition $(S_{0})$. We have
$|\psi''(s)| <l$ putting $N=\max |\psi'(s)|$ and choosing $|\alpha |
\leq \alpha_{0}= \min \left\{ 1, \frac{k}{2(N+l)},
\frac{\rho}{N+l^{2}}\right\}$. Then we have, for $|\alpha| \leq
\alpha_{0}$, that
\[
t_{\alpha}'(s)= t'(s)+\alpha \psi'(s) +\alpha^{2} \psi''(s) \geq
k-(N+l)\alpha_{0} \geq k-\frac{k}{2}= \frac{k}{2}> 0 \, ,
\]
$s \in E_{1}\bigcup E_{2}$, and $C_{\alpha}$ has an absolutely
continuous representation $x=x_{\alpha}(t), a \leq t \leq b$. We
also have  $|t_{\alpha}(s)-t(s)|<|\alpha|(N+l^{2}) <  \rho$. Hence
$|x_{\alpha}(t)-x(t)|= |x(t_{\alpha}(s))-x(t(s))|< \delta$ and we
conclude that $J[C_{\alpha}]\geq J[C]$.  On the other hand, by
setting $\beta(\alpha, s) = F(t, X, t', X', t'', X'')$, we have by
differentiation that
$$ \left.\frac{\partial \beta}{\partial
\alpha}\right|_{\alpha = 0}= \frac{\partial F}{\partial t} \psi +
\frac{\partial F}{\partial t'} \psi'+  \frac{\partial F}{\partial
t''} \psi'' \, ,
$$
where
\begin{equation}
\label{eq5}
\begin{split}
\frac{\partial F}{\partial t} &= \frac{\partial \mathcal{L}}{\partial t}t'  \, , \\
\frac{\partial F}{\partial t'}&= \mathcal{L}- \frac{\partial
\mathcal{L}}{\partial \dot{x}}\frac{\dot{x}}{t'}+\frac{1}{t'^{2}}
\frac{\partial \mathcal{L}}{\partial
\ddot{x}}\left(\frac{-2\ddot{x}}{t'}
+\frac{3t''}{t'}\dot{x}\right) \, , \\
\frac{\partial F}{\partial t''}&= -\frac{1}{t'^{2}}\frac{\partial
\mathcal{L}}{\partial \ddot{x}}\dot{x} \, .
\end{split}
\end{equation}
By hypotheses $(S_{i})_{0\leq i \leq n}$  both terms $\frac{\partial
F}{\partial t}\psi, \frac{\partial F}{\partial t'}\psi',
\frac{\partial F}{\partial t''}\psi''$  are bounded  in
$E_{1}\bigcup E_{2}$ by a function which is $L$-integrable in $[0,l]$.
Then, we can differentiate under the sign of
the integral to obtain:
$$
0= \left.\frac{\partial J(C_{\alpha})}{d \alpha}\right|_{\alpha= 0}
= \int_{0}^{l}\left( \frac{\partial F}{\partial t} \psi +
\frac{\partial F}{\partial t'} \psi'+  \frac{\partial F}{\partial
t''} \psi'' \right)ds \, .
$$
Integration by parts, and using \eqref{eq3}--\eqref{eq4}, yields
\begin{equation*}
\begin{split}
0 &= \int_{0}^{l} \phi_{0}(s)\psi'' ds =
\int_{E_{1}}\phi_{0}(s)\psi'' ds
+ \int_{E_{2}}\phi_{0}(s)\psi'' ds \\
& \leq |E_{1}| |E_{2}|(d_{1}-d_{2}) < 0
\end{split}
\end{equation*}
which is a contradiction. Equality \eqref{eq2} is now proved.
\end{proof}

%%%%%%%%%%%%%%%%%%%%%%%%%%%%%%%%%%%%%%%%%%%%

\subsection{Generalized Euler-Lagrange equations}

Arguments similar to those used to prove Theorem~\ref{theorem1}
can be utilized to prove a generalized Euler-Lagrange equation. This
condition is not necessary in the proof of our regularity theorem, but
is given here because of its significance: necessary conditions for
\eqref{P2} in the class $W^{2}_{2}$ have an interest of their own
(\textrm{cf.} Example~\ref{ex:CV90}).

\begin{theorem}
\label{theorem2} Under the hypotheses $(S_{i})_{1\leq i \leq n}$,
if $x(\cdot) \in W_2^2$ is a minimizer of problem \eqref{P2}, then
we have the following integral form of the Euler-Lagrange equations:
\begin{equation}
\phi_{i}(s)= \frac{\partial F}{\partial \ddot{x}_{i}}(\theta(s))
- \int_{0}^{s}\frac{\partial F}{\partial \dot{x}_{i}}(\theta(\sigma)) d\sigma
+ \int_{0}^{s}\int_{0}^{\tau}\frac{\partial F}{\partial {x}_{i}}(\theta(\sigma)) d\sigma d\tau
= c_{i}, \quad 1\leq i \leq n
 , \label{eq6}
\end{equation}
where functions $\frac{\partial F}{\partial \ddot{x}_{i}},
\frac{\partial F}{\partial \dot{x}_{i}}, \frac{\partial F}{\partial
{x}_{i}}$ are evaluated at $\theta(s) = \left(t(s), X(s), t'(s), X'(s), t''(s),
X''(s)\right)$ and $c_{i}$, $i \in \{1,\ldots,n\}$, denote constants.
\end{theorem}

\begin{proof}
The proof is also by contradiction and is similar to that of
Theorem\ref{theorem1}. Suppose that \eqref{eq6} is not satisfied.
For $i=1 \ldots n$ and $|\alpha| \leq 1$, we consider the curve
$C_{\alpha}: t= t_{\alpha}(s)$, $x= X_{\alpha}(s)$, $0 \leq s \leq
l$, with
\begin{gather*}
X_{i \alpha}(s) = X_{i}(s) + \alpha \psi(s) +\alpha^{2} \psi'(s) \, ,\\
X_{j \alpha}(s)= X_{j}(s), \quad j \neq i \, .
\end{gather*}
We have $|\psi''(s)| \leq l$ $a.e$ and, if we put $N=
\max|\psi'(s)|$, then for $$|\alpha| \leq \alpha_{0}= \min \left\{ 1,
\frac{\delta}{(N+1)l}, \frac{\delta}{N+l}, \frac{\delta}{l}\right\}$$
we can write that
$$
\left|X_{i \alpha}(s)-X_{i}(s)\right| =\left|\alpha \psi +\alpha^{2}
\psi' \right| \leq \alpha_{0}(N+1)l \leq \delta \, ,
$$
$$
\left|\dot{X}_{i \alpha}(s)-\dot{X_{i}}(s)\right| =\left|\alpha
\psi' +\alpha^{2} \psi'' \right| \leq \alpha_{0}(N+l) \leq \delta \,
,
$$
$$
\left|\ddot{X}_{i \alpha}(s)-\ddot{X}_{i}(s)\right| = \left|\alpha
\psi'' \right| \leq \alpha_{0}l \leq \delta  \, ,
$$
and thus $J[C] \leq J[C_{\alpha}]$ for all  $|\alpha|\leq
\alpha_{0}$. Setting, as before,
$$\beta(\alpha, s)=  F(t(s), X(s), t'(s), X'(s), t''(s), X''(s))
$$
we have
$$ \left.\frac{\partial \beta}{\partial
\alpha}\right|_{\alpha = 0}= \frac{\partial F}{\partial x_{i}} \psi
+ \frac{\partial F}{\partial \dot{x}_{i}} \psi'+  \frac{\partial
F}{\partial \ddot{x}_{i}} \psi'', \mbox{ for } s\in [0, l] \quad
 a. e.
$$
Note that by the hypotheses $(S_{i})_{1 \leq i \leq n}$
$$
\left|\frac{\partial F}{\partial x_{i}} \right| =
\left|\frac{\partial \mathcal{L}}{\partial x_{i}}t' \right| \leq
G\left(t(s), X, \frac{X'}{t'},
\frac{1}{t'^{2}}(X''-\frac{X'}{t'})t''\right)t' \, ,
$$
$$
\left|\frac{\partial F}{\partial \dot{x}_{i}}\right|
 \leq
G\left(t(s),X, \frac{X'}{t'},
\frac{1}{t'^{2}}(X''-\frac{X'}{t'})t''\right) + G\left(t, X,
\frac{X'}{t'}, \frac{1}{t'^{2}}(X''-\frac{X'}{t'})t'' \right)
\frac{t''}{t'^{2}} \, ,
$$
$$
\left|\frac{\partial F}{\partial \ddot{x}_{i}} \right|=\left|
\frac{\partial \mathcal{L}}{\partial \ddot{x}_{i}} \right|\leq
G\left(t(s),X, \frac{X'}{t'},
\frac{1}{t'^{2}}(X''-\frac{X'}{t'})t''\right) \frac{1}{t'}
$$
are $L$-integrable in $[0, l]$. Thus, the terms $\frac{\partial
F}{\partial x_{i}} \psi , \frac{\partial F}{\partial \dot{x}_{i}}
\psi', \frac{\partial F}{\partial \ddot{x}_{i}} \psi''$ are bounded
in $E_{1}\bigcup E_{2}$ by a fixed $L$-integrable function.
For $s \in (E_{1}\bigcup E_{2})^{c}$, we
have $\psi''(s)=0$ and $\frac{\partial F}{\partial \ddot{x}_{i}}
\psi''= 0$. The proof is continued in the same lines as in the end of
the proof of Theorem~\ref{theorem1}, applying the usual
rule of differentiation under the integral sign and integration by
parts, which leads to a contradiction.
\end{proof}

%%%%%%%%%%%%%%%%%%%%%

\section{Regularity result for autonomous problems}
\label{sec:Reg}

We shall present now a regularity result for \eqref{P2} under
certain additional requirements on the Lagrangian $\mathcal{L}$.

\begin{theorem}
\label{theorem4} In addition to the hypotheses $(S_{i})_{0\leq i
\leq n}$, let us consider the autonomous problem \eqref{P2},
\textrm{i.e.} let us assume that $\mathcal{L}$ does not depend on
$t$: $\mathcal{L} = \mathcal{L}(x,\dot{x},\ddot{x})$. If
$\frac{\partial \mathcal{L}}{\partial \ddot{x}}$ is superlinear,
\textrm{i.e.} there exist constants $a >0$ and $b > 0$ such that
\begin{equation}
\label{H4} a |w|+ b \leq \left|\frac{\partial \mathcal{L}}{\partial
\ddot{x}}(s, v, w)\right| \mbox{ for all } (s, v, w)\in
\mathbb{R}^{n}\times \mathbb{R}^{n} \times \mathbb{R}^{n}  \, ,
\end{equation}
then every minimizer $x \in W^{2}_{2}$ of the problem is on
$W_{2}^{\infty}$.
\end{theorem}
\begin{remark} Theorem~\ref{theorem4} remain valid if instead
of the superlinearity condition \eqref{H4} we impose the stronger
quadratically coercive condition:
there exist constants $a >0$ and $b > 0$ such that
\begin{equation*}
 a |w|^{2}+ b \leq \left|\frac{\partial \mathcal{L}}{\partial
\ddot{x}}(s, v, w)\right| \mbox{ for all } (s, v, w)\in
\mathbb{R}^{n}\times \mathbb{R}^{n} \times \mathbb{R}^{n} \, .
\end{equation*}
\end{remark}

\begin{example}
A trivial example of a Lagrangian satisfying all the conditions
$(S_{i})_{0\leq i \leq n}$ and \eqref{H4} is $\mathcal{L}(x,
\dot{x}, \ddot{x})= \mathcal{L}(\ddot{x}) = a \ddot{x}^2 + b
\ddot{x}$ with $a$ and $b$ strictly positive constants (one can
choose $G(t, x, \dot{x}, \ddot{x}) = 2 a |\ddot{x}| + b$ $\in L^{2}$
in $(S_{i})$). It follows from Theorem~\ref{theorem4} that all
minimizers of the problem
\begin{gather*}
I[x(\cdot)] = \int_{t_0}^{t_1} \left[a \ddot{x}(t)^2 + b \ddot{x}(t)\right] dt \longrightarrow \min \\
x(\cdot) \in W^{2}_{2} \, , \quad a , b > 0 \\
x(t_0) = \alpha \, , \quad x(t_1) = \beta
\end{gather*}
are $W_{2}^{\infty}$ functions.
\end{example}

As an immediate corollary to our Theorem~\ref{theorem4}, we obtain
conditions of non-occurrence of the Lavrentiev phenomenon for the
autonomous second-order variational problems.
\begin{corollary}
Under the hypotheses of Theorem~\ref{theorem4}, the autonomous
problems do not admit the Lavrentiev gap $W^{2}_{2}-W_{2}^{\infty}$:
\begin{equation*}
\inf_{x(\cdot) \in W^{2}_{2}} \int_{a}^{b} \mathcal{L}\left(x(t),\dot{x}(t) ,
\ddot{x}(t)\right) dt = \inf_{x(\cdot) \in W_{2}^{\infty}}
\int_{a}^{b} \mathcal{L}\left(x(t),\dot{x}(t) , \ddot{x}(t)\right) dt \, .
\end{equation*}
\end{corollary}

\begin{example}
\label{ex:CV90} Let us consider the autonomous problem proposed in
\cite{cla,clavin90} ($n=1$, $m = 2$): $\mathcal{L}(s, v, w)=
\left|s^{2}-v^{5}\right|^{2}|w|^{22}+\varepsilon |w|^{2}$, $t \in [0,1]$. The
problem satisfies hypotheses (H1)-(H3) of Tonelli's existence
theorem. Function $\tilde{x}(t)= k t^{\frac{5}{3}}$ verifies the
integral form of the Euler-Lagrange equations \eqref{eq6}. However,
$\tilde{x}$ belongs to $W^{2}_{2}$ but not to $W_{2}^{\infty}$. The
regularity condition \eqref{H4} of Theorem~\ref{theorem4} is not
satisfied.
\end{example}

\begin{proof}(of Theorem~\ref{theorem4})
Using  \eqref{eq2} and \eqref{eq5} we get
\begin{equation*}
\frac{1}{t'^{2}} \frac{\partial \mathcal{L}}{\partial \ddot{x}}
\dot{x} + \int_{0}^{s}\left \{ \mathcal{L}- \frac{1}{t'}
\frac{\partial \mathcal{L}}{\partial \dot{x}}
\dot{x}+\frac{1}{t'^{2}} \frac{\partial \mathcal{L}}{\partial
\ddot{x}} \left(\frac{-2\ddot{x}}{t'}+\frac{3t''}{t'}\dot{x} \right) \right \}
- \int\int \frac{\partial \mathcal{L}}{\partial t}t'= c_{0}
\end{equation*}
and since we are in the autonomous case,
\begin{equation*}
\frac{1}{t'^{2}} \frac{\partial \mathcal{L}}{\partial \ddot{x}}
\dot{x} + \int_{0}^{s}\left \{ \mathcal{L}- \frac{1}{t'}
\frac{\partial \mathcal{L}}{\partial \dot{x}}
\dot{x}+\frac{1}{t'^{2}} \frac{\partial \mathcal{L}}{\partial
\ddot{x}} \left(\frac{-2\ddot{x}}{t'}+\frac{3t''}{t'}\dot{x} \right) \right \}
= c_{0} \, .
\end{equation*}
Therefore,
\begin{equation*}
\begin{split}
\frac{1}{t'^{2}} \frac{\partial \mathcal{L}}{\partial \ddot{x}}
\dot{x} &= c_{0}- \int_{0}^{s}\left \{ \mathcal{L}- \frac{1}{t'}
\frac{\partial \mathcal{L}}{\partial \dot{x}}
\dot{x}+\frac{1}{t'^{2}} \frac{\partial \mathcal{L}}{\partial
\ddot{x}}
\left(\frac{-2\ddot{x}}{t'}+\frac{3t''}{t'}\dot{x} \right) \right \}\\
&= c_{0}- \int_{0}^{s} \mathcal{L} + \int_{0}^{s}\frac{1}{t'}
\frac{\partial \mathcal{L}}{\partial \dot{x}} \dot{x}+ 2
\int_{0}^{s} \frac{1}{t'^{3}}\frac{\partial \mathcal{L}}{\partial
\ddot{x}} \ddot{x} - \int_{0}^{s} \frac{3t''}{t'^{3}} \frac{\partial
\mathcal{L}}{\partial \ddot{x}}\dot{x} \, .
\end{split}
\end{equation*}
Applying the Holder's inequality, we obtain
\begin{equation*}
\left|\frac{1}{t'^{2}} \frac{\partial \mathcal{L}}{\partial
\ddot{x}}\dot{x} \right| \leq |c_{0}|+ \| \mathcal{L}\|_{1}
 + k_{1}\left\| \frac{\partial \mathcal{L}}{\partial \dot{x}}\right\|_{2} \left\|\dot{x}\right\|_{2}
 +  k_{2}\left\|\frac{\partial \mathcal{L}}{\partial \ddot{x}}\right\|_{2}\|\ddot{x}\|_{2}
 + \int_{0}^{s} \left|\frac{3t''}{t'}\right| \,
 \left|\frac{1}{t'^{2}} \frac{\partial \mathcal{L}}{\partial \ddot{x}}
\dot{x} \right| \, ,
\end{equation*}
where $k_{1}, k_{2}$ are  positive constants. Then, using the fact
that $\mathcal{L}\in C^{1}, \mathcal{L}, \frac{\partial
\mathcal{L}}{\partial \dot{x}}, \frac{\partial \mathcal{L}}{\partial
\ddot{x}} \in L^{2}$ and $x \in W_{2}^{2}$ (in other terms, $x,
\dot{x}, \ddot{x}\in L^{2}$), it follows that $\frac{1}{t'^{2}}
\frac{\partial \mathcal{L}}{\partial \ddot{x}}\dot{x}$ satisfies a
condition of the form
$$
\left|\frac{1}{t'^{2}} \frac{\partial \mathcal{L}}{\partial
\ddot{x}}\dot{x} \right| \leq k_{3}+ \int_{0}^{s} \left|\frac{3t''}{t'}\right|
\left|\frac{1}{t'^{2}} \frac{\partial \mathcal{L}}{\partial
\ddot{x}} \dot{x} \right| \, ,
$$
for a certain  positive constant $k_{3}$. Now, Gronwall's Lemma
leads to the  following uniform bound:
$$ \left|\frac{1}{t'^{2}} \frac{\partial \mathcal{L}}{\partial \ddot{x}}\dot{x} \right| \leq  k_{4} $$
with a positive constant $k_{4}$. Since $ t' \leq 1$, we deduce that
$\frac{\partial \mathcal{L}}{\partial \ddot{x}}\dot{x} $  is
uniformly bounded. Besides, since $\left|\frac{\partial
\mathcal{L}}{\partial \ddot{x}}\right|$ verifies \eqref{H4}, we have
$$
|\dot{x}| \left(a |\ddot{x}|+ b\right) \leq \left| \frac{\partial
\mathcal{L}}{\partial \ddot{x}}\dot{x} \right| \leq k_{4} \quad
(b>0) \, .
$$
Therefore we get for a positive constant $k_{5}$
$$
|\dot{x}| \leq \frac{k_{4}}{a |\ddot{x}|+ b}  \leq k_{5} \, .
$$ Then $\frac{\partial \mathcal{L}}{\partial \ddot{x}}$ is
uniformly bounded.
Since $\frac{\partial \mathcal{L}}{\partial \ddot{x}}(s, v, w)$
goes to $+ \infty$ with $|w|$ (by superlinearity), this implies a uniform bound
 on $|\ddot{x}|$ which leads to the intended conclusion
that $\ddot{x}$ is essentially bounded.
\end{proof}

Theorems~\ref{theorem1}, \ref{theorem2} and \ref{theorem4} admit a
generalization for problems of an order higher than two. This is under
study and will be addressed in a forthcoming paper.

%%%%%%%%%%%%%%%%%%%%%

\section{Conclusions}

The search for appropriate conditions on the data of the problems of
the calculus of variations with higher-order derivatives, under
which we have regularity of solutions or under which more general
necessary conditions hold, is an important area of study. In this
paper we have obtained necessary optimality conditions of
duBois-Reymond and Euler-Lagrange type, valid in the class of
functions where the existence is proved. Minimizers in this class
may have unbounded derivatives and fail to satisfy the classical
necessary conditions of duBois-Reymond or Euler-Lagrange. We prove
that if the derivatives of the Lagrangian function with respect to
the highest derivatives verify a superlinear condition, then all the
minimizers have essentially bounded derivatives. This imply
non-occurrence of the Lavrentiev phenomenon and validity of
classical necessary optimality conditions.

%%%%%%%%%%%%%%%%%%%%%

\section*{Acknowledgements}

The authors are grateful to Ilona Dzenite and Enrique H. Manfredini
for the suggestions regarding improvement of the text;
to Andrei V. Sarychev for pointing out several mistakes in an earlier
version of the manuscript. Possible remaining errors are, of course, the
sole responsibility of the authors. The first author was supported by
the \emph{Portuguese Foundation for Science and Technology} (FCT),
through project SFRH/BPD/20934/2004; the second author
by the \emph{Centre for Research in Optimization and Control}
(CEOC) from FCT, cofinanced by the European Community Fund FEDER/POCTI.

%%%%%%%%%%%%%%%%%%%%%%%%%%%%%%%%%%%%%%%%%%%%%%%%%%%%%%%

\end{document}